\documentclass[12pt]{amsart}

\usepackage{calc}  
\usepackage{amsmath}  
\usepackage{amsthm}
\usepackage{amsfonts}
\usepackage{amssymb}
\usepackage{amstext}
\usepackage{amscd}
\usepackage{graphicx}
\usepackage{enumerate}
\usepackage[all]{xy}

\newtheorem{theorem}{Theorem}[section]
\newtheorem{proposition}[theorem]{Proposition}
\newtheorem{corollary}[theorem]{Corollary}
\newtheorem{lemma}[theorem]{Lemma}

\theoremstyle{definition}

\newtheorem{definition}[theorem]{Definition}
\newtheorem{remark}[theorem]{Remark}

\newtheorem{acknow}{Acknowledgment}

\def\N{\mathbb{N}}
\def\Z{\mathbb{Z}}
\def\Q{\mathbb{Q}}
\def\R{\mathbb{R}}
\def\C{\mathbb{C}}

\def\bE{\mathbb{E}}
\def\bP{\mathbb{P}}

\def\CP{\mathbb{C} \mathbb{P}}
\def\barCP2{\overline{\mathbb{CP}}\ \! \!^2}

\def\<{\left\langle}
\def\>{\right\rangle}
\def\({\left(}
\def\){\right)}

\def\cA{\mathcal{A}}
\def\cB{\mathcal{B}}

\def\cG{\mathcal{G}}
\def\cH{\mathcal{H}}

\def\cL{\mathcal{L}}

\def\cM{\mathcal{M}}

\def\cV{\mathcal{V}}

\def\fg{\mathfrak{g}}

\def\id{\operatorname{id}}

\def\Tr{\operatorname{Tr}}

\begin{document}

\title{An SO(3)-version of 2-torsion instanton invariants}
\author[H. Sasahira]{Hirofumi Sasahira$^*$}
\date{}

\renewcommand{\thefootnote}{\fnsymbol{footnote}}
\footnote[0]{2000\textit{ Mathematics Subject Classification}. 57R57}

\thanks{
$^*$Partially supported by the 21th century COE program
at Graduate School of Mathematical Sciences,
the University of Tokyo.
}

\address{
Graduate School of Mathematical Sciences,
University of Tokyo,\endgraf 
3-8-1 Komaba Meguro-ku, Tokyo 153-8941, Japan
}
\email{sasahira@ms.u-tokyo.ac.jp}

\maketitle


\begin{abstract}
We construct an invariant for non-spin $4$-manifolds by using $2$-torsion cohomology classes of moduli spaces of instantons on $SO(3)$-bundles. The invariant is an $SO(3)$-version of Fintushel-Stern's $2$-torsion instanton invariant. We show that this $SO(3)$-torsion invariant is non-trivial for $2\CP^2 \# \barCP2$, while it is known that any known invariant of $2\CP^2 \# \barCP2$ coming from the Seiberg-Witten theory is trivial since $2\CP^2 \# \barCP2$ has a positive scalar curvature metric.
\end{abstract}

\section{Introduction}
The purpose of this paper is to construct an $SO(3)$-version of Fintushel-Stern's torsion invariants \cite{FS}. R. Fintushel and R. Stern constructed a variant of Donaldson invariants for spin $4$-manifolds by using $2$-torsion cohomology classes of the moduli spaces of instantons on $SU(2)$-bundles.  They used cohomology classes of degree one and two.  S. K. Donaldson gave another construction by using a class of degree $3$ \cite{yang}.
As is well known, the usual Donaldson invariant is trivial for the connected sum of $4$-manifolds with $b^+$ positive (\cite{poly}). On the other hand, Fintushel and Stern showed that their torsion invariant is not necessarily trivial for the connected sum of the form $Y \# S^2 \times S^2$ in general. 

In this paper, we define an invariant of $4$-manifolds using $2$-torsion cohomology classes of $SO(3)$-moduli spaces and show that our invariant is not necessarily trivial for $Y \# S^2 \times S^2$ as in the case of Fintushel-Stern's invariant. We basically follow the argument in \cite{FS} and modify it to extend the definition to non-spin $4$-manifolds.

The outline of the construction is as follows. Let $X$ be a closed, oriented, simply connected, non-spin Riemannian $4$-manifold and $P$ be an $SO(3)$-bundle over $X$ satisfying
\[
w_2(P)=w_2(X) \in H^2(X; \Z_2), \quad p_1(P) \equiv \sigma(X) \mod 8.
\]
Here $\sigma(X)$ is the signature of $X$. Let $\cB_P^*$ be the space of gauge equivalence classes of irreducible connections on $P$. In \cite{AMR}, S. Akbulut, T. Mrowka and Y. Ruan showed that $H^1(\cB_P^*; \Z_2)$ is isomorphic to $\Z_2$. We denote the generator by $u_1$. On the other hand, for homology class $[\Sigma] \in H_2(X; \Z)$ with self-intersection number even, we have an integral cohomology class $\mu([\Sigma]) \in H^2(\cB_P^*; \Z)$. Suppose that the dimension of the moduli space $M_P$ of instantons on $P$ is $2d+1$ for some non-negative integer $d$. In general $M_P$ is not compact. However for homology classes $[\Sigma_1], \dots, [\Sigma_d] \in H_2(X; \Z)$ with self-intersection numbers even, we can define the pairing
\[
q_{X}^{u_1}([\Sigma_1], \dots, [\Sigma_d])=
\< u_1 \cup \mu([\Sigma_1]) \cup \cdots \cup \mu([\Sigma_d]), [M_P] \> \in \Z_2
\]
in an appropriate sense. We show that this number depends only on the homology classes $[\Sigma_i]$ and gives a differential-topological invariant of $X$. 

We will show a gluing formula of torsion invariants for $Y \# S^2 \times S^2$, which is an $SO(3)$-version of Theorem 1.1 in \cite{FS}. By using this gluing formula and D. Kotschick's calculation in \cite{K, K2}, we prove that $q_{2\CP^2 \# \barCP2}^{u_1}$ is non-trivial. This example exhibits two interesting aspects explained below.

The first aspect is related to vanishing theorem.
We have a description of $X=2\CP^2 \# \barCP2$ as the connected sum of $Y_1=\CP^2$ and $Y_2=\CP^2 \# \barCP2$. Since the second Stiefel-Whitney class $w_2(P)$ is equal to $w_2(X)$, both of $w_2(P)|_{Y_1}$ and $w_2(P)|_{Y_2}$ are non-trivial. In such a situation, the usual Donaldson invariants are trivial by the dimension-count argument (\cite{MM}). Hence the non-triviality of $q_{2\CP^2 \# \barCP2}^{u_1}$ implies that the dimension-count argument can not be applied directly to proving such a vanishing theorem in our case. If each homology class $[\Sigma_i]$ is in $H_2(Y_1; \Z)$ or $H_2(Y_2; \Z)$, then we can show that our invariant vanishes. However we can not reduce the argument to this case because of the condition that $[\Sigma_i] \cdot [\Sigma_i]$ must be even to define our invariant.

The next aspect is related to the Seiberg-Witten theory. In \cite{Wi}, E. Witten introduced invariants, called the Seiberg-Witten invariants, of $4$-manifolds using monopole equations. He conjectured that the invariants are equivalent to the Donaldson invariants and explicitly wrote a formula which should give a relation between the Donaldson invariants and the Seiberg-Witten invariants. In \cite{PT}, V. Pidstrigach and A. Tyurin proposed a program to give a rigorous mathematical proof of the formula by using non-abelian monopoles. The theory of non-abelian monopoles has been developed by P. Feehan and T. Leness (\cite{FL1, FL2, FL3}). Feehan and Leness recently announced that they completed the proof of Witten's formula for $4$-manifolds of simple type with $b_1=0$ and $b^+ > 1$ in \cite{FL4}.

The non-triviality of $q^{u_1}_{2\CP^2 \# \barCP2}$ is quite a contrast to the equivalence of the Donaldson invariants and Seiberg-Witten invariants. 
If a $4$-manifold has a positive scalar curvature metric and satisfies $b^+(X) \geq 1$, then the moduli space of solutions of the monopole equations with respect to the metric is empty for some perturbation. Hence any known invariant of $2\CP^2 \# \barCP2$ coming from the monopole equations (the Seiberg-Witten invariant and a refinement due to S. Bauer and M. Furuta \cite{BF}) is trivial since $2\CP^2 \# \barCP2$ has a positive scalar curvature metric.

The paper is organized as follows. In Section 2, we construct cohomology classes $\mu([\Sigma])$ and $u_1$, and define a torsion invariant. In Section 3, we prove a gluing formula for the connected sum of the form $Y \# S^2 \times S^2$. In Section 4, we prove that $q_{2 \CP^2 \# \barCP2}^{u_1}$ is non-trivial by using the gluing formula. We also discuss the reason why the usual vanishing theorem does not hold for our torsion invariant.

\begin{acknow}
The author would like to thank my advisor Mikio Furuta for his suggestions and warm encouragement. The author also thanks Yukio Kametani and Nobuhiro Nakamura for useful conversations.
\end{acknow}

\section{Torsion invariants} \label{invariants}

\subsection{Notations} \label{notations}
Let $X$ be a closed, oriented, simply connected $4$-manifold, $g$ a Riemannian metric on $X$ and $P$ an $SO(3)$-bundle over $X$. Put
\[
k=-\frac{1}{4}p_1(P) \in \Q, \quad w=w_2(P) \in H^2(X; \Z_2).
\]
Let $\cA_P^*$ be the space of irreducible connections on $P$ and $\cG_P$ be the gauge group of $P$. We write $\cB_P^*$ or $\cB_{k,w,X}^*$ for the quotient space $\cA_P^*/\cG_P$. We denote by $M_P$ or $M_{k,w,X}$  the moduli space of instantons on $P$.

Let $A$ be an instanton on $P$. We have a sequence
\[
\Omega_X^0(\fg_{P}) \stackrel{d_A}{\longrightarrow} \Omega^1_X (\fg_P) \stackrel{d_A^+}{\longrightarrow} \Omega_X^+(\fg_P).
\]
The condition that $A$ is an instanton implies that $d_A^+ \circ d_A = 0$. Hence the above sequence define a complex. We denote the cohomology groups by $H_{A}^0$, $H_A^1$, $H_A^2$.

Let $\bar{P}$ be a $U(2)$-lift of $P$ and $\bar{E}$ be  the rank $2$ complex vector bundle associated with $\bar{P}$. Fix a connection $a_{\det}$ on $\det \bar{E}$. We write $\cA_{\bar{E}}$ for the space of connections on $\bar{E}$ which induce the connection $a_{\det}$ on $\det \bar{E}$, and write $\cA_{\bar{E}}^*$ for the space of irreducible connections in $\cA_{\bar{E}}$. Let $\cG_{\bar{E}}$ be the group of bundle automorphisms on $\bar{E}$ with determinant $1$. We also introduce a subgroup $\cG_{\bar{E}}^0$ of $\cG_{\bar{E}}$. Fix a point $x_0$ in $X$. The subgroup $\cG_{\bar{E}}^0$ is defined by
\[
\cG_{\bar{E}}^0=\{ g \in \cG_{\bar{E}} | g(x_0)=1 \}.
\]
We denote the quotient spaces by
\[
\cB_{\bar{E}}^* = \cA_{\bar{E}}^* / \cG_{\bar{E}}, \quad
\widetilde{\cB}_{\bar{E}} = \cA_{\bar{E}} / \cG_{\bar{E}}^0, \quad
\widetilde{\cB}_{\bar{E}}^* = \cA_{\bar{E}}^* / \cG_{\bar{E}}^0.
\]
Since we are assuming that $X$ is simply connected, the natural map $\cB^*_{\bar{E}} \rightarrow \cB^*_{P}$ is bijective.

To construct cohomology classes $u_1$ and $\mu([\Sigma])$, we need the universal bundle $\widetilde{\bE}$ over $X \times \widetilde{\cB}_{\bar{E}}$. The universal bundle is defined by
\[
\widetilde{\bE}:=\bar{E} \times_{\cG_{\bar{E}}^0} \cA_{\bar{E}} \longrightarrow X \times \widetilde{\cB}_{\bar{E}}.
\]
For a closed, oriented surface $\Sigma$ embedded in $X$, let $\nu(\Sigma)$ be a small tubular neighborhood of $\Sigma$. We define spaces of gauge equivalence classes of connections on $\nu(\Sigma)$.  Let $\cA_{\nu(\Sigma)}$ be the space of connections on $\bar{E}|_{\nu(\Sigma)}$ which induce the connection $a_{\det}|_{\nu(\Sigma)}$ on $\det \bar{E}|_{\nu(\Sigma)}$. Let $\cG_{\nu(\Sigma)}$ be the group of automorphisms of $\bar{E}|_{\nu(\Sigma)}$ with determinant $1$. We assume that the  base point $x_0$ is in $\nu(\Sigma)$. Define  $\cG_{\nu(\Sigma)}^0$ by
\[
\cG_{\nu(\Sigma)}^0 = \{ g \in \cG_{\nu(\Sigma)} | g(x_0)=1 \}.
\]
We denote the quotient spaces by
\[
\cB_{\nu(\Sigma)}^* = \cA_{\nu(\Sigma)}^* / \cG_{\nu(\Sigma)}, \quad
\widetilde{\cB}_{\nu(\Sigma)} = \cA_{\nu(\Sigma)} / \cG_{\nu(\Sigma)}^0, \quad
\widetilde{\cB}_{\nu(\Sigma)}^* = \cA_{\nu(\Sigma)}^* / \cG_{\nu(\Sigma)}^0.
\]
Restricting connections, we have a map
\[
\tilde{r}_{\nu(\Sigma)}:\widetilde{\cB}_{\bar{E}}^* \longrightarrow \widetilde{\cB}_{\nu(\Sigma)}.
\]
We have the universal bundle $\widetilde{\bE}_{\nu(\Sigma)}$ over $\nu(\Sigma) \times \widetilde{\cB}_{\nu(\Sigma)}$ defined by
\[
\widetilde{\bE}_{\nu(\Sigma)}:=(\bar{E}|_{\nu(\Sigma)}) \times_{\cG_{\nu(\Sigma)}^0} \cA_{\nu(\Sigma)} \longrightarrow \nu(\Sigma) \times \widetilde{\cB}_{\nu(\Sigma)}.
\]

\subsection{Cohomology classes of $\cB_{P}^*$} \label{cohomolgy}
Suppose $\Sigma$ is a closed, oriented surface embedded in $X$ such that $\< w_2(P),[\Sigma] \> \equiv 0 \mod 2$.  In this subsection, we define a $2$-dimensional integral cohomology class $\mu([\Sigma]) \in H^2(\cB_P^*; \Z)$. Basically we follow a standard construction in \cite{DK, K}.

We first define the cohomology class $\tilde{\mu}_{\bar{E}}([\Sigma]) \in H^2(\widetilde{\cB}_{\bar{E}}; \Z)$ to be the slant product $c_2(\widetilde{\bE})/[\Sigma]$.

\begin{lemma} \label{lem beta}
Let $\beta:\widetilde{\cB}_{\bar{E}}^* \rightarrow \cB_{\bar{E}}^*$ be the projection.  Then the induced homomorphism
\[
\beta^*:H^2(\cB_{\bar{E}}^*; \Z) \longrightarrow H^2( \widetilde{\cB}_{\bar{E}}^*; \Z)
\]
is injective.
Moreover for a homology class $[\Sigma] \in H_2(X; \Z)$ with $\< w_2(P), [\Sigma] \> \equiv 0 \mod 2$, 
the cohomology class $\tilde{\mu}_{\bar{E}}([\Sigma])$ lies in the image of $\beta^*$.
\end{lemma}

\begin{proof}
Since $H^1(SO(3); \Z)=0$, the spectral sequence associated with the fibration $SO(3) \rightarrow \widetilde{\cB}_{\bar{E}}^* \rightarrow \cB_{\bar{E}}^*$ induces an exact sequence
\begin{equation} \label{eq exact}
0 \longrightarrow H^2(\cB_{\bar{E}}^*; \Z) \stackrel{\beta^*}{\longrightarrow} H^2(\widetilde{\cB}_{\bar{E}}^*; \Z) \longrightarrow H^2(SO(3); \Z),
\end{equation}
which implies the injectivity of $\beta^*$. 

Let $\eta$ be a complex line bundle over $SO(3)$ defined by
\[
\eta:=SU(2) \times_{ \{ \pm 1 \} } \C \longrightarrow SO(3).
\]
Here the action of $\{ \pm 1 \}$ on $\C$ is the multiplication. Then it is easy to obtain the identification
\[
\widetilde{\bE}|_{\Sigma \times SO(3)} = (\bar{E}|_{\Sigma}) \times_{ \{ \pm 1 \} } SU(2) = (\bar{E}|_{\Sigma}) \boxtimes \eta \longrightarrow \Sigma \times SO(3),
\]
and we have
\[
\begin{split}
c_2(\widetilde{\bE}|_{\Sigma \times SO(3)})/[\Sigma]
&=\big( \pi_1^* c_2(\bar{E}|_{\Sigma}) + \pi_1^* c_1(\bar{E}|_{\Sigma}) \cup \pi_2^* c_1(\eta) \big)/[\Sigma] \\
&=\< c_1(\bar{E}), [\Sigma] \> c_1(\eta) \\
&\in H^2(SO(3); \Z) \cong \Z_2,
\end{split}
\]
where 
\[
\pi_1:\Sigma \times SO(3) \longrightarrow \Sigma, \quad
\pi_2:\Sigma \times SO(3) \longrightarrow SO(3)
\]
are the projections.
If $\< w_2(P), [\Sigma] \>$ is zero, the pairing $\< c_1(\bar{E}), [\Sigma] \>$ is even, and hence the restriction of $c_2(\widetilde{\bE})/[\Sigma]$ to $SO(3)$ is trivial. From the exact sequence (\ref{eq exact}), $\tilde{\mu}_{\bar{E}}([\Sigma])$ is in the image of $\beta^*$.
\end{proof}

By Lemma \ref{lem beta}, there is a unique element of $H^2(\cB_{\bar{E}}^*; \Z)$ such that the image by $\beta^*$ is $\tilde{\mu}_{\bar{E}}([\Sigma])$. Through the natural identification between $\cB_{P}^*$ and $\cB_{\bar{E}}^*$, we have a $2$-dimensional cohomology class of $\cB_{P}^*$. We denote it by $\mu_{\bar{E}}([\Sigma])$.

\begin{lemma} \label{lem mu}
Let $X$ be a closed, oriented, simply connected $4$-manifold and $P$ be an $SO(3)$-bundle over $X$. Suppose that $[\Sigma]$ is a $2$-dimensional homology class in $X$ with $\< w_2(P),[\Sigma] \> \equiv 0 \mod 2$.
Then the cohomology class $\mu_{\bar{E}}([\Sigma]) \in H^2(\cB_P^*; \Z)$ is independent of the choice of $\bar{E}$.
\end{lemma}

This lemma will be shown in \S \ref{well-def} as a corollary of Lemma \ref{lem pull back}.
Under the assumption in Lemma \ref{lem mu}, we define $\mu([\Sigma]) \in H^2(\cB_P^*; \Z)$ as follows.

\begin{definition} \label{def mu}
For a homology class $[\Sigma] \in H_2(X, \Z)$ with $\< w_2(P), [\Sigma] \> \equiv 0 \mod 2$, the cohomology class $\mu([\Sigma]) \in H^2(\cB_P^*; \Z)$ is defined to be $\mu_{\bar{E}}([\Sigma])$.
\end{definition}

\begin{remark}
Let
\[
\bP:=P \times_{\cG_P} \cA^*_P  \longrightarrow X \times \cB_P^*
\]
be the universal bundle of $P$. Then the usual definition of $\mu$-map is given by
\[
\begin{array}{rccc}
\mu_{\Q}: & H_2(X; \Z) & \longrightarrow & H^2(\cB_P^*; \Q) \\
     & [\Sigma]  & \longmapsto & -\frac{1}{4}p_1(\bP)/[\Sigma].
\end{array}
\]
In general, $\mu_{\Q}([\Sigma])$ does not have an integral lift. Under our assumptions, it is easy to see that $\mu([\Sigma])$ is an integral lift of $\mu_{\Q}([\Sigma])$.
\end{remark}

Next we define a torsion cohomology class $u_1 \in H^1(\cB_P^*; \Z_2)$. We write $\sigma(X)$ for the signature of $X$. Akbulut, Mrowka and Ruan showed the following in \cite{AMR}.

\begin{proposition} [\cite{AMR}] \label{pi1}
Let $X$ be a closed, oriented, simply connected $4$-manifold and $P$ be an $SO(3)$-bundle over $X$. Then we have
\[
\pi_1(\cB_P^*)=\left\{
\begin{array}{cl}
\Z_2 &  \text{if $w_2(P)=w_2(X),\ p_1(P) \equiv \sigma(X) \mod 8$} \\
  1  & \text{otherwise}.
\end{array}
\right.
\]
\end{proposition}

\begin{remark} \label{rem pi1}
Suppose  $P$ is an $SO(3)$-bundle over $X$ with $w_2(P)$ equal to $w_2(X)$ and let $\bar{P}$ be a $U(2)$-lift of $P$. Then $p_1(P)$ is equal to $\sigma(X)$ modulo $8$ if and only if $c_2(\bar{P})$ is equal to $0$ modulo $2$. This equivalence is a consequence of the formulas
\[
p_1(P)=-4c_2(\bar{P}) + c_1(\bar{P})^2, \quad w_2(X)^2 \equiv \sigma(X) \mod 8.
\]
\end{remark}

When $w_2(P)=w_2(X)$ and $p_1(P) \equiv \sigma(X) \mod 8$, we have $H^1(\cB_P^*; \Z_2) \cong \Z_2$ from Proposition \ref{pi1}.

\begin{definition}
Let $X$ be a closed, oriented, simply connected $4$-manifold and $P$ be an $SO(3)$-bundle over $X$ satisfying
$w_2(P) = w_2(X), \quad p_1(P) \equiv \sigma(X) \mod 8$.
We write $u_1$ for the generator of $H^1(\cB_P^*; \Z_2) \cong \Z_2$.
\end{definition}

\subsection{Construction of $q_{X}^{u_1}$}
Let $X$ be a closed, oriented, simply connected $4$-manifold. Suppose $b^+(X)=2a$ for a positive integer $a$. Let $P$ be an $SO(3)$-bundle over $X$.
Assume that $P$ satisfies the condition
\begin{equation} \label{eq w_2}
w_2(P)=w_2(X) \in H^2(X; \Z_2), \quad p_1(P) \equiv \sigma(X) \mod 8. 
\end{equation}
The virtual dimension of $M_P$ is given by
\[
\dim M_P=-2p_1(P)-3(1+b^+(X))=8k-3(1+2a).
\]
If we put $d=-p_1(P)-3a-2=4k-3a-2$, then we have
\[
\dim M_P=2d+1.
\]
From the condition (\ref{eq w_2}), we have
\[
d \equiv -\sigma(X) -3a - 2 \mod 8.
\]
Suppose that $d \geq 0$ and take $2$-dimensional homology classes $[\Sigma_1], \dots, [\Sigma_d]$ of $X$ satisfying
\[
\< w_2(P), [\Sigma_i] \> \equiv 0 \mod 2 \quad (i=1,\dots, d).
\]
The assumption $\< w_2(P),[\Sigma_i] \> \equiv 0 \mod 2$ is equivalent to  $[\Sigma_i] \cdot [\Sigma_i] \equiv 0 \mod 2$ since $w_2(P)$ is equal to $w_2(X)$. We want to define the pairing $\< u_1 \cup \mu([\Sigma_1]) \cup \cdots \cup \mu([\Sigma_d]), M_P \> \in \Z_2$. The moduli space $M_P$ is not compact in general and the pairing is not well-defined in the usual sense.
To define the pairing, we need submanifolds $V_{\Sigma_i}$ dual to $\mu([\Sigma_i])$ which behave nicely near the ends of $M_P$. We briefly explain how the submanifolds are constructed. See \cite{poly, DK} for the details. 

We use the following three things.
The first is that when $b^+(X)$ and $k=-\frac{1}{4}p_1(P)$ are positive $M_P$ lies in $\cB_P^*$ and has a natural smooth structure for generic metrics on $X$.
The second is that the restrictions of irreducible instantons to open subsets are also irreducible.
The third is that the cohomology class $\mu([\Sigma])$ comes from $\cB_{\nu(\Sigma)}^*$. More precise statement of the third is as follows.

Let $[\Sigma] \in H_2(X; \Z)$ be a homology class with $[\Sigma] \cdot [\Sigma] \equiv 0 \mod 2$.
Since the following diagram is commutative
\[
\begin{CD}
\widetilde{\bE}|_{\nu(\Sigma) \times \cB_{\bar{E}}^*}=(\bar{E}|_{\nu(\Sigma)}) \times_{\cG_{\bar{E}}^0} \cA^*_{\bar{E}} 
@>{\id_{\bar{E}} \times \tilde{r}_{\nu(\Sigma)}}>> 
\widetilde{\bE}_{\nu(\Sigma)}=(\bar{E}|_{\nu(\Sigma)}) \times_{\cG_{\nu(\Sigma)}^0} \cA_{\nu(\Sigma)} \\
@VVV @VVV \\
\nu(\Sigma) \times \widetilde{\cB}_{\bar{E}}^* 
@>>{\id_{\nu(\Sigma)} \times \tilde{r}_{\nu(\Sigma)}}> 
\nu(\Sigma) \times \widetilde{\cB}_{\nu(\Sigma)}   
\end{CD}
\]
we obtain
\begin{equation} \label{eq mu restriction}
\tilde{\mu}_{\bar{E}}([\Sigma])=c_2(\tilde{\bE})/[\Sigma]=\tilde{r}_{ \nu(\Sigma)}^* (c_2(\widetilde{\bE}_{\nu(\Sigma)})/[\Sigma]) \in H^2(\widetilde{\cB}_{\bar{E}}^*; \Z).
\end{equation}
We apply Lemma \ref{lem beta} to the restriction $P|_{\nu(\Sigma)}$, instead of $P$ itself. Then we see that there exists a unique $2$-dimensional cohomology class $\mu_{\nu(\Sigma),\bar{E}}([\Sigma])$ of $\cB_{\nu(\Sigma)}^*$ such that the pull-back by the natural projection $\widetilde{\cB}_{\nu(\Sigma)}^* \rightarrow \cB_{\nu(\Sigma)}^*$ is equal to $c_2(\widetilde{\bE}_{\nu(\Sigma)})/[\Sigma]$. 

We define $V_{\Sigma}$ as follows.

\begin{definition}
Take a homology class $[\Sigma] \in H_2(X; \Z)$ with $[\Sigma] \cdot [\Sigma]$ even.
We write $\cL_{\Sigma}$ for a complex line bundle over $\cB_{\nu(\Sigma), \bar{E}}^*$ with first Chern class $\mu_{\nu(\Sigma), \bar{E}}([\Sigma]) \in H^2(\cB_{\nu(\Sigma),\bar{E}}^*; \Z)$. Fix a section $s_{\Sigma}$ of $\cL_{\Sigma}$. We denote the zero locus of $s_{\Sigma}$ by $V_{\Sigma} \subset \cB_{\nu(\Sigma)}^*$. Suppose that $b^+(X)$ and $k=-\frac{1}{4}p_1(P)$ are positive. For a generic metric $g$, we define
\[
M_{P} \cap V_{\Sigma} := \{ \ [A] \in M_P \ | \ [A|_{\nu(\Sigma)}] \in V_{\Sigma} \ \}.
\]
\end{definition}
We will show that the pairing $\< u_1, M_P \cap V_{\Sigma_1} \cap \cdots \cap V_{\Sigma_d} \>$ is well-defined under some condition.

\begin{remark} \label{rem line bundle}
We give some remarks on the line bundle $\cL_{\Sigma}$. We refer to \cite{poly, DK} for details.
\begin{itemize}
\item
As is well-known, we are also able to construct the line bundle $\cL_{\Sigma}$ by using a family of twisted Dirac operators on $\Sigma$.

\item
Assume that $\< w_2(P), [\Sigma] \>$ is equal to $0$ modulo $2$. Then $P|_{\nu(\Sigma)}$ is topologically trivial.
Let $\cB_{\nu(\Sigma) \ +}^* := \cB_{\nu(\Sigma)}^* \cup \{ [\Theta_{\nu(\Sigma)} ] \}$. Here $\Theta_{\nu(\Sigma)}$ is the trivial connection on $\nu(\Sigma)$. It is known that $\cL_{\Sigma}$ extends to $\cB_{\nu(\Sigma) \ +}^*$.  Hence we  can assume that the section $s_{\Sigma}$ is non-zero near $[\Theta_{\nu(\Sigma)}]$.
In the case when $w_2(P)$ is zero, we need this property to define invariants. On the other hand, when we treat an $SO(3)$-bundle $P$ with $w_2(P)$ non-trivial, we do not need this property for the definition of invariants. However we will need this property in Lemma \ref{lem bubble} to prove some property of our invariant .
\end{itemize}
\end{remark}

We prepare some lemmas. The following is well-known.

\begin{lemma} [\cite{conn, DK}] \label{lem transverse} 
Let $X$ be a closed, oriented, simply connected $4$-manifold with $b^+(X)$ positive and $P$ be an $SO(3)$-bundle with $w_2(P)=w_2(X)$ and $k=-\frac{1}{4}p_1(P)$ positive. Take homology classes $[\Sigma_1], \dots, [\Sigma_{d'}] \in H_2(X; \Z)$ with self-intersection numbers even. For generic sections $s_{\Sigma_i}$, the intersections
\[
M_{k-j,w,X} \cap \left( \bigcap_{i \in I} V_{\Sigma_i} \right)
\quad
(I \subset \{ 1, \dots , d' \},\ 0 \leq j < k)
\]
are transverse.
\end{lemma}

From now on, we require that $\Sigma_i$ are generic in the following sense.

\begin{equation} \label{eq transverse}
\left\{ \quad
\begin{split}
\Sigma_i \pitchfork \Sigma_j \quad & \text{($i,j$ distinct)}  \\
\Sigma_i \cap \Sigma_j \cap \Sigma_k=\emptyset \quad  & \text{($i,j,k$ distinct)}.
\end{split}
\right.
\end{equation}

\begin{lemma} \label{lem compact}
Let $X$ be a closed, oriented, simply connected, non-spin $4$-manifold with $b^+(X)$ positive. Let $P$ be an $SO(3)$-bundle over $X$ with $w_2(P)$ equal to $w_2(X)$. Suppose that the dimension of $M_P$ is $2d'+r$ for a non-negative integer $d'$ and $1 \leq r \leq 3$. 
Take $d'$ homology classes $[\Sigma_1],\dots, [\Sigma_{d'}] \in H^2(X; \Z)$ with
\[
[\Sigma_i] \cdot [\Sigma_i] \equiv 0 \mod 2 \quad (i=1,\dots,d').
\]
Moreover we assume that the surfaces $\Sigma_i$ satisfy the condition (\ref{eq transverse}).
Then for generic sections $s_{\Sigma_i}$, the intersection
\[
M_P \cap V_{\Sigma_1} \cap \cdots \cap V_{\Sigma_{d'}}
\]
is a compact $r$-dimensional manifold.
\end{lemma}

\begin{proof}
Put $k=-\frac{1}{4}p_1(P)$, $w=w_2(P)$. For $[A] \in M_P$, we have
\[
\begin{split}
k &= -\frac{1}{4}p_1(P) \\
&=\frac{1}{8\pi^2} \int_X \Tr(F_A^2) \\
&=\frac{1}{8\pi^2} \int_X | F_{A}^- |^2 d\mu_g  - \frac{1}{8 \pi^2} \int_X | F_{A}^+ |^2 d\mu_g \\
&=\frac{1}{8\pi^2} \int_X | F_A^- | d\mu_g \geq 0.
\end{split}
\]
by the Chern-Weil theory.
Here $d\mu_g$ is the volume form with respect to $g$. First we show $k>0$. If not, $k=0$ and $A$ is flat. Since $X$ is simply connected, $A$ is trivial. This contradicts to the assumption that $w_2(P)$ is non-trivial. Hence we have $k>0$. From Lemma \ref{lem transverse},  $M_P \cap V_{\Sigma_1} \cap \cdots \cap V_{\Sigma_{d'}}$ is a smooth $r$-dimensional manifold for generic sections $s_{\Sigma_i}$.

Next we prove that $M_P \cap V_{\Sigma_1} \cap \cdots \cap V_{\Sigma_{d'}}$ is compact.
Let $\{ [A^{(n)}] \}_{n \in \N}$ be a sequence in $M_P \cap V_{\Sigma_1} \cap \cdots \cap V_{\Sigma_{d'}}$. Uhlenbeck's weak compactness theorem implies that there is a subsequence $\{ [A^{(n')}] \}_{n'}$ which is weakly convergent to
\[
([A_{\infty}]; x_1,\dots, x_l) \in M_{k-l, w, X} \times X^l.
\]
We also have $k-l>0$ in the same way as above. Let $m$ be the number of the tubular neighborhoods $\nu(\Sigma_i)$ which contain $x_{\alpha}$ for some $\alpha$ with $1 \leq \alpha \leq l$. Then without loss of generality, we may suppose that
\[
[A_{\infty}] \in M_{k-l, w, X} \cap V_{\Sigma_1} \cap \cdots \cap V_{\Sigma_{d'-m}}
\]
if we change the order of the surfaces. If we take the tubular neighborhoods $\nu(\Sigma_i)$ to be sufficiently small, we have
\begin{equation*} \label{eq nu sigma}
\nu(\Sigma_i) \cap \nu(\Sigma_j) \cap \nu(\Sigma_k) = \emptyset
\quad
(\text{$i, j, k$ distinct})
\end{equation*}
from (\ref{eq transverse}). Hence we have $m \leq 2l$.
Since $k-l>0$, the intersection $M_{k-l,x,X} \cap V_{\Sigma_1} \cap \cdots \cap V_{\Sigma_{d'-m}}$ is transverse by Lemma \ref{lem transverse}. From this transversality, we obtain
\[
\begin{split}
0 
&\leq \dim M_{k-l, w, X} \cap V_{\Sigma_1} \cap \cdots \cap V_{\Sigma_{d'-m}} \\
&= \dim M_{k, w, X} - 8l - 2(d'-m) \\
&=r-8l+2m \\
&\leq r-4l.
\end{split}
\]
Since we suppose $1 \leq r \leq 3$, we have $l=0$ and
\[
[A_{\infty}] \in M_{P} \cap V_{\Sigma_1} \cap \cdots \cap V_{\Sigma_{d'}}.
\]
Hence $M_P \cap V_{\Sigma_1} \cap \cdots \cap V_{\Sigma_{d'}}$ is compact.
\end{proof}

Let $X$ be as in Lemma \ref{lem compact} and $P$ be an $SO(3)$-bundle over $X$ satisfying (\ref{eq w_2}). Suppose that $\dim M_P$ is $2d+1$ for a non-negative integer $d$ and take homology classes $[\Sigma_1], \dots, [\Sigma_d] \in H_2(X; \Z)$ with self-intersection numbers even. From Lemma \ref{lem compact}, we have the pairing
\[
\< u_1, M_P \cap V_{\Sigma_1} \cap \cdots \cap V_{\Sigma_d} \> \in \Z_2.
\]

\begin{proposition} \label{prop well-defined}
Let $X$ be a closed, oriented, simply connected, non-spin $4$-manifold with $b^+(X)=2a$ for a positive integer $a$ and $P$ be an $SO(3)$-bundle over $X$ satisfying (\ref{eq w_2}). Assume that the dimension of $M_P$ is $2d+1$ for a non-negative integer $d$. Then the pairing 
\begin{equation*} \label{eq pairing}
\< u_1, M_P \cap V_{\Sigma_1} \cap \cdots \cap V_{\Sigma_d} \> \in \Z_2
\end{equation*}
is independent of the choices of Riemannian metric $g$, $U(2)$-lift $\bar{P}$ of $P$, sections $s_{\Sigma_i}$ of $\cL_{\Sigma_i}$ and surfaces $\Sigma_i$ representing the homology classes $[\Sigma_i]$. Moreover the pairing is multi-linear with respect to $[\Sigma_1], \dots, [\Sigma_d]$.
\end{proposition}

We prove the above proposition in \S \ref{well-def}. By using this proposition, we can easily show that the following invariant $q^{u_1}_X$ is well defined.

\begin{definition}
Let $X$ be as in Proposition \ref{prop well-defined}.
Let $A_d'(X)$  be the subspace of $\otimes^d H^2(X;\Z)$ generated by 
\[
\{ \ [\Sigma_1] \otimes \cdots \otimes [\Sigma_d] \ | \ [\Sigma_i] \in H_2(X; \Z), [\Sigma_i] \cdot [\Sigma_i] \equiv 0 \mod 2 \ \},
\]
and we put
\[
A'(X) := \bigoplus_{d} A'_d(X),
\] 
where $d$ runs over non-negative integers with $d \equiv -\sigma(X)-3a-2 \mod 8$. We define $q_X^{u_1}$ by
\[
\begin{array}{rccc}
q_X^{u_1}:& A'(X) & \longrightarrow & \Z_2 \\
          & ([\Sigma_1],\dots,[\Sigma_d]) & \longmapsto & q_{k,w,X}^{u_1}([\Sigma_1],\dots,[\Sigma_d]) := \< u_1,M_P \cap V_{\Sigma_1} \cap \cdots \cap V_{\Sigma_d} \>.
\end{array}
\]
Here $P$ is an $SO(3)$-bundle over $X$ with $w_2(P)=w_2(X)$ and $p_1(P)=-d-3a-2$.
\end{definition}

\subsection{Well-definedness of $q_{X}^{u_1}$} \label{well-def}

In this subsection, we prove Proposition \ref{prop well-defined}.
First we show the independence of $q_X^{u_1}$ from Riemannian metric $g$ and sections $s_{\Sigma_i}$ in a standard way. Take two metrics $g$, $g'$ on $X$ and sections $s_{\Sigma_i}$, $s_{\Sigma_i}'$ of $\cL_{\Sigma_i}$. Choose a path $\{ g_t \}_{t \in [0,1]}$ between $g$ and $g'$, and a path $\{ s_{\Sigma_i,t} \}_{t \in [0,1] }$ between $s_{\Sigma_i}$ and $s_{\Sigma_i}'$. Then put
\[
\cM:= \coprod_{t \in [0,1]} M_{P}(g_t) \times \{ t \}, \quad
\cM \cap \cV_{\Sigma_i} := \{ \ ([A], t) \in \cM \ | \ s_{\Sigma_i, t}([A|_{\nu(\Sigma_i)}]) = 0 \ \}.
\]
Using a similar argument in the proof of Lemma \ref{lem compact}, we can show the following lemma:

\begin{lemma} \label{lem wd}
Let $X$ and $P$ be as in Proposition \ref{prop well-defined}. Then for generic paths $\{ g_{t} \}_{t \in [0,1]}$ and $\{ s_{\Sigma_i,t} \}_{t \in [0,1]}$, the intersection
\[
\cM \cap \cV_{\Sigma_1} \cap \cdots \cap \cV_{\Sigma_d}
\]
is a compact $2$-dimensional manifold whose boundary is
\[
( M_{P}(g) \cap V_{\Sigma_1} \cap \cdots \cap V_{\Sigma_d} ) \coprod
( M_{P}(g') \cap V_{\Sigma_1}' \cap \cdots \cap V_{\Sigma_d}').
\]
\end{lemma}
This lemma implies
\[
\< u_1, M_P(g) \cap V_{\Sigma_1} \cap \cdots \cap V_{\Sigma_d} \>=
\< u_1, M_P(g') \cap V_{\Sigma_1}' \cap \cdots \cap V_{\Sigma_d}' \>
\in \Z_2,
\]
and the pairing $\< u_1, M_P \cap V_{\Sigma_1} \cap \cdots \cap V_{\Sigma_d} \>$ is independent of the choices of $g$ and $s_{\Sigma_i}$.

Next we see the independence of $q_{X}^{u_1}$ from the choice of $U(2)$-lift $\bar{P}$ of $P$ .
Take two $U(2)$-lifts $\bar{P}$ and $\bar{P}'$ of $P$. The associated vector bundle $\bar{E}'$ with $\bar{P}'$ is topologically isomorphic to $\bar{E} \otimes L$ for some complex line bundle $L$  over $X$. 
Fix connections $a_{\det}$, $a_L$ on $\det \bar{E}$, $L$ and an isomorphism
\[
\varphi:\bar{E}'  \stackrel{\cong}{\longrightarrow} \bar{E} \otimes L.
\]
We have a connection  $a_{\det}'$ on $\det \bar{E}'$ induced by $a_{\det}, a_L$ and $\varphi$. We consider connections on $\bar{E} \otimes L$ and $\bar{E}'$ which are compatible with $a_{\det}+2a_L$ and $a_{\det}'$ respectively. By tensoring $a_L|_{\nu(\Sigma)}$, we have maps
\[
t_{\cA}:\cA_{\nu(\Sigma), \bar{E}} \stackrel{\cong}{\longrightarrow} \cA_{\nu(\Sigma), \bar{E} \otimes L}, \quad
t_{\cB^*}:\cB_{\nu(\Sigma), \bar{E}}^* \stackrel{\cong}{\longrightarrow} \cB_{\nu(\Sigma), \bar{E} \otimes L}^*, \quad
t_{\widetilde{\cB}}:\widetilde{\cB}_{\nu(\Sigma), \bar{E}} \stackrel{\cong}{\longrightarrow} \widetilde{\cB}_{\nu(\Sigma), \bar{E} \otimes L}.
\]
Moreover the pull-back by $\varphi$ induces identifications
\[
\psi_{\cB^*}:\cB_{\nu(\Sigma), \bar{E} \otimes L}^* \stackrel{\cong}{\longrightarrow} \cB_{\nu(\Sigma), \bar{E}'}^*, \quad
\psi_{\widetilde{\cB}}:\widetilde{\cB}_{ \nu(\Sigma), \bar{E} \otimes L} \stackrel{\cong}{\longrightarrow} \widetilde{\cB}_{\nu(\Sigma), \bar{E}'}.
\]

\begin{lemma} \label{lem pull back}
Suppose $\cL_{\Sigma}$, $\cL_{\Sigma}'$ are complex line bundles over $\cB_{\nu(\Sigma), \bar{E}}^*$, $\cB_{\nu(\Sigma), \bar{E}'}^*$ corresponding to the cohomology classes $\mu_{\nu(\Sigma), \bar{E}}([\Sigma]) \in H^2(\cB_{\nu(\Sigma),\bar{E}}^*; \Z)$, $\mu_{\nu(\Sigma), \bar{E}'}([\Sigma]) \in H^2(\cB_{\nu(\Sigma),\bar{E}'}^*; \Z)$. Then we have
\[
(\psi_{\cB^*} \circ t_{{\cB}^*})^* \cL_{\Sigma}' \cong \cL_{\Sigma}.
\]
\end{lemma}

\begin{proof}
It is sufficient to show that $(\psi_{\widetilde{\cB}} \circ t_{\widetilde{\cB}})^* (c_2(\widetilde{\bE}'_{\nu(\Sigma)})/[\Sigma])$ is equal to $c_2(\widetilde{\bE}_{\nu(\Sigma)})/[\Sigma]$ since $H^2(\cB^*_{\nu(\Sigma),\bar{E}}; \Z) \rightarrow H^2(\widetilde{\cB}^*_{\nu(\Sigma), \bar{E}}; \Z)$ is injective.

Let $\pi_1:\nu(\Sigma) \times \widetilde{\cB}_{\nu(\Sigma), \bar{E}} \rightarrow \nu(\Sigma)$ be the projection. 
We have the following commutative diagram:
\[
\begin{CD}
\widetilde{\bE}_{\nu(\Sigma)} \otimes \pi_1^* (L|_{\nu(\Sigma)})@.
\widetilde{\bE}_{\nu(\Sigma)}' \\
@| @| \\
(\bar{E} \otimes L|_{\nu(\Sigma)}) \times_{\cG^0_{\nu(\Sigma),\bar{E}}} \cA_{\nu(\Sigma),\bar{E}} 
@>{\varphi^{-1} \times ( \varphi^* \ \circ \ t_{\cA} )}>> 
(\bar{E}'|_{\nu(\Sigma)}) \times_{\cG^0_{\nu(\Sigma),\bar{E}'}} \cA_{\nu(\Sigma),\bar{E}'} \\
@VVV      @VVV \\
\nu(\Sigma) \times \widetilde{\cB}_{\nu(\Sigma),\bar{E}}
@>>{\id_{\nu(\Sigma)} \times ( \psi_{\widetilde{\cB}} \ \circ \ t_{\tilde{\cB}} )}>
\nu(\Sigma) \times \widetilde{\cB}_{\nu(\Sigma),\bar{E}'}
\end{CD}
\]
Hence we have
\[
\big( \id_{\nu(\Sigma)} \times (\psi_{\widetilde{\cB}} \circ t_{\tilde{\cB}}) \big)^* \ \widetilde{\bE}'_{\nu(\Sigma)} \cong 
\widetilde{\bE}_{\nu(\Sigma)} \otimes \pi_1^* (L|_{\nu(\Sigma)})
\]
and we obtain
\begin{equation} \label{eq slant}
\begin{split}
(\psi_{\widetilde{\cB}} \circ t_{\widetilde{\cB}})^* (c_2(\widetilde{\bE}'_{\nu(\Sigma)})/[\Sigma])
&=c_2(\widetilde{\bE}_{\nu(\Sigma)} \otimes \pi_1^* (L|_{\nu(\Sigma)}))/[\Sigma] \\
&=\big\{ c_2(\widetilde{\bE}_{\nu(\Sigma)}) + \pi_1^* c_1(L|_{\nu(\Sigma)}) \cup c_1(\widetilde{\bE}_{\nu(\Sigma)}) + \pi_1^* c_1(L|_{\nu(\Sigma)})^2 \big\} / [\Sigma] \\
&=c_2(\widetilde{\bE}_{\nu(\Sigma)})/[\Sigma] + \big\{  \pi_1^* c_1(L|_{\nu(\Sigma)}) \cup c_1(\widetilde{\bE}_{\nu(\Sigma)}) \big\} /[\Sigma] \\
&\in H^2(\widetilde{\cB}_{\bar{E}}; \Z).
\end{split}
\end{equation}
By the K\"unneth formula, we can write
\[
\begin{split}
c_1(\widetilde{\bE}_{\nu(\Sigma)}) 
&=c_1(\widetilde{\bE}_{\nu(\Sigma)})_{\nu(\Sigma)} + 
c_1(\widetilde{\bE}_{\nu(\Sigma)})_{\widetilde{\cB}} \\
&\in H^2(\nu(\Sigma) \times \widetilde{\cB}_{\nu(\Sigma), \bar{E}}; \Z) \cong
H^2(\nu(\Sigma); \Z) \oplus H^2(\widetilde{\cB}_{\nu(\Sigma), \bar{E}}; \Z)
\end{split}
\]
since $\widetilde{\cB}_{\nu(\Sigma), \bar{E}}$ is simply connected (\cite{AB}). The action of $\cG_{\nu(\Sigma), \bar{E}}^0$ on $\Lambda^2 \bar{E}|_{\nu(\Sigma)}$ is trivial, since the determinants of elements of $\cG^0_{\nu(\Sigma), \bar{E}}$ are equal to $1$ by definition. Hence $\Lambda^2 \widetilde{\bE}_{\nu(\Sigma)}$ is the pull-back $\pi_1^* (\Lambda^2 \bar{E}|_{\nu(\Sigma)})$. This implies that the $\widetilde{\cB}_{\nu(\Sigma)}$-part $c_1(\widetilde{\bE}_{\nu(\Sigma)})_{\widetilde{\cB}}$ of $c_1(\widetilde{\bE}_{\nu(\Sigma)})=c_1(\Lambda^2 \widetilde{\bE}_{\nu(\Sigma)})$ is $0$ and we have
\[
\begin{split}
\big\{  \pi_1^* c_1(L|_{\nu(\Sigma)}) \cup c_1(\widetilde{\bE}_{\nu(\Sigma)}) \big\} /[\Sigma]
&=\big\{ \pi_1^* c_1(L|_{\nu(\Sigma)}) \cup c_1(\widetilde{\bE}_{\nu(\Sigma)})_{\nu(\Sigma)} \big\}  / [\Sigma]\\
&=0 \in H^2(\widetilde{\cB}_{\nu(\Sigma)}; \Z).
\end{split}
\]
From the equation (\ref{eq slant}), we obtain
\begin{equation} \label{eq c2}
(\psi_{\widetilde{\cB}} \circ t_{\widetilde{\cB}})^* (c_2(\widetilde{\bE}'_{\nu(\Sigma)})/[\Sigma]) =
c_2(\widetilde{\bE}_{\nu(\Sigma)})/[\Sigma] 
\in H^2(\widetilde{\cB}_{\nu(\Sigma)}; \Z).
\end{equation}
\end{proof}

\noindent
{\it Proof of Lemma \ref{lem mu}.} \\
Lemma \ref{lem mu} follows from (\ref{eq c2}) and the following commutative diagram:
\[
\xymatrix{
\widetilde{\cB}_{\nu(\Sigma),\bar{E}} \ar[rr]^{\cong}_{\psi_{\widetilde{\cB}} \circ t_{\widetilde{\cB}}} & &
\widetilde{\cB}_{\nu(\Sigma), \bar{E}'} \\
\widetilde{\cB}^*_{X,\bar{E}} \ar[u]^{\tilde{r}_{\nu(\Sigma)}} \ar[rr]^{\cong} \ar[d] &  &
\widetilde{\cB}^*_{X, \bar{E}'} \ar[u]_{\tilde{r}_{\nu(\Sigma)}} \ar[d] \\
\cB_{X,\bar{E}}^* \ar[rr]^{\cong} \ar[dr]^{\cong} & &
\cB_{X,\bar{E}'}^* \ar[dl]_{\cong} \\
& \cB_{X,P}^*
}
\]

\qed

\noindent
{\it Proof of independence of $q_{X}^{u_1}$ from $\bar{P}$.} \\
Take homology classes $[\Sigma_i] \in H_2(X; \Z)$ with $[\Sigma_i] \cdot [\Sigma_i] \equiv 0 \mod 2$ for $i=1,\dots,d$ and choose $U(2)$-lifts $\bar{P}$ and $\bar{P}'$ of $P$. Then we obtain line bundles $\cL_{\Sigma_i}$ and $\cL_{\Sigma_i}'$ over $\cB_{\nu(\Sigma_i), \bar{E}}^*$ and $\cB_{\nu(\Sigma_i), \bar{E}'}^*$.
We denote the zero locus of  sections $s_{\Sigma_i}$, $s_{\Sigma_i}'$ of $\cL_{\Sigma_i}$, $\cL_{\Sigma_i}'$ by $V_{\Sigma_i}$, $V_{\Sigma_i}'$.
By Lemma \ref{lem pull back}, $(\psi_{\cB^*} \circ t_{{\cB}^*})^* \cL_{\Sigma_i}'$ is isomorphic to $\cL_{\Sigma_i}$.  We fix an isomorphism and regard the section $s_{\Sigma_i}'$ of $\cL_{\Sigma_i}'$ as a sections of $\cL_{\Sigma_i}$ through the identifications 
\[
\psi_{\cB^*} \circ t_{{\cB}^*}:\cB_{\nu(\Sigma_i), \bar{E}}^* \stackrel{\cong}{\longrightarrow} \cB^*_{\nu(\Sigma_i), \bar{E}'}, \quad
(\psi_{\cB^*} \circ t_{\cB^*})^* \cL_{\Sigma_i}' \cong \cL_{\Sigma_i}. 
\]
We take paths $\{ s_{\Sigma_i,t} \}_{t \in [0,1]}$ between $s_{\Sigma_i}$ and $s_{\Sigma_i}'$. In the same way as Lemma \ref{lem wd}, we have a bordism between $M_P \cap V_{\Sigma_1} \cap \cdots \cap V_{\Sigma_d}$ and $M_P \cap V_{\Sigma_1}' \cap \cdots \cap V_{\Sigma_d}'$. Hence we obtain
\[
\< u_1, M_P \cap V_{\Sigma_1} \cap \cdots \cap V_{\Sigma_d} \> =
\< u_1, M_P \cap V_{\Sigma_1}' \cap \cdots \cap V_{\Sigma_d}' \> \in \Z_2.
\]
\qed

Lastly we show that $q_{X}^{u_1}$ is independent of the choice of surfaces $\Sigma_i$ representing the homology classes $[\Sigma_i]$ and that $q_{X}^{u_1}$ is multi-linear with respect to $[\Sigma_1], \dots, [\Sigma_d]$. It follows from the following lemma directly.

\begin{lemma}
Let $X$ and $P$ be as in Proposition \ref{prop well-defined}.
Take homology classes $[\Sigma_1], \dots, [\Sigma_d] \in H_2(X;\Z)$  with self-intersection numbers even. Moreover assume that
\[
[\Sigma_1] = [\Sigma_1'] + [\Sigma_1''] \in H_2(X;\Z), \quad
[\Sigma_1'] \cdot [\Sigma_1'] \equiv [\Sigma_1''] \cdot [\Sigma_1''] \equiv 0 \mod 2.
\]
Then we have
\begin{gather*}
\< u_1,M_P \cap V_{\Sigma_1} \cap V_{\Sigma_2} \cap \cdots \cap V_{\Sigma_d} \>=\\
\< u_1,M_P \cap V_{\Sigma_1'} \cap V_{\Sigma_2} \cap \cdots \cap V_{\Sigma_d} \>+ \< u_1,M_P \cap V_{\Sigma_1''} \cap V_{\Sigma_2} \cap \cdots \cap V_{\Sigma_d} \>
\in \Z_2.
\end{gather*}
\end{lemma}

\begin{proof}
By definition, we have
\[
\tilde{\mu}_{\bar{E}}([\Sigma_1])
=c_2(\tilde{\bE})/[\Sigma_1]
=c_2(\tilde{\bE})/[\Sigma_1'] + c_2(\tilde{\bE})/[\Sigma_1'']
=\tilde{\mu}_{\bar{E}}([\Sigma_1'])+\tilde{\mu}_{\bar{E}}([\Sigma_1''])
\in H^2(\widetilde{\cB}_{\bar{E}}; \Z).
\]
The homomorphism $\beta^*:H^2(\cB_{\bar{E}}^*; \Z) \rightarrow H^2(\widetilde{\cB}_{\bar{E}}^*; \Z)$ is injective and $\tilde{\mu}_{\bar{E}}([\Sigma_1])$, $\tilde{\mu}_{\bar{E}}([\Sigma_1'])$, $\tilde{\mu}_{\bar{E}}([\Sigma_1''])$ lie in the image $\beta^*$ from Lemma \ref{lem beta}. Hence we have
\[
\mu([\Sigma_1])=\mu([\Sigma_1'])+\mu([\Sigma_1'']) \in H^2(\cB_P^*; \Z).
\]
Since $M_{P} \cap V_{\Sigma_2} \cap \cdots \cap V_{\Sigma_d}$ is compact from Lemma \ref{lem compact}, we have
\[
\begin{split}
&\< u_1, M_P \cap V_{\Sigma_1} \cap \cdots \cap V_{\Sigma_d} \> \\
&\quad = \< u_1 \cup \mu([\Sigma_1]), M_P \cap V_{\Sigma_2} \cap \cdots \cap V_{\Sigma_d} \> \\
&\quad =\< u_1 \cup (\mu([\Sigma_1']) + \mu([\Sigma_1''])), M_P \cap V_{\Sigma_2} \cap \cdots \cap V_{\Sigma_d} \> \\
&\quad =\< u_1 \cup \mu([\Sigma_1']), M_P \cap V_{\Sigma_2} \cap \cdots \cap V_{\Sigma_d} \> + \< u_1 \cup \mu([\Sigma_1'']), M_P \cap V_{\Sigma_2} \cap \cdots \cap V_{\Sigma_d} \> \\
&\quad =\< u_1, M_P \cap V_{\Sigma_1'} \cap V_{\Sigma_2} \cap \cdots \cap V_{\Sigma_d} \> + \< u_1, M_P \cap V_{\Sigma_1''} \cap V_{\Sigma_2} \cap \cdots \cap V_{\Sigma_d} \>.
\end{split}
\]
\end{proof}

\section{A connected sum formula for $Y \# S^2 \times S^2$}
\subsection{Statement of the result}
As is well known Donaldson invariants vanish for the connected sum $X_1 \# X_2$ provided $b^+(X_i)>0$ for $i=1, 2$ (\cite{poly}). In \cite{FS}, however, Fintushel and Stern defined  some torsion invariants by using instantons on $SU(2)$-bundles and they showed that their $SU(2)$-torsion invariants are non-trivial for the connected sum of the form $Y \# S^2 \times S^2$. In this section, we show a similar non-vanishing theorem for our $SO(3)$-torsion invariants.

Let $Y$ be a closed, oriented, simply connected, non-spin $4$-manifold with $b^+(Y)=2a-1$ for $a>1$. Let $Q$ be an $SO(3)$-bundle with $w_2(Q)$ equal to $w_2(Y)$ and $p_1(Q)$ equal to $\sigma(Y)+4$ modulo $8$. Suppose that the dimension of $M_Q$ is $2d$ for a non-negative integer $d$. When we fix an orientation on the space $\cH_g^+(Y)$ of self-dual harmonic $2$-forms on $Y$ and an lift $c \in H^2(Y;\Z)$ of $w_2(Q) \in H^2(Y;\Z_2)$,  we have the Donaldson invariant
\[
q_{k-1,w,Y}:\otimes^d H_2(Y;\Z) \longrightarrow \Q
\]
where
\[
k-1=-\frac{1}{4}p_1(Q) \in \Q, \quad w=w_2(Q) \in H^2(Y;\Z_2).
\]
When $[\Sigma_i] \cdot [\Sigma_i]$ are even for $i=1,\dots, d$, then $q_{k-1,w,Y}([\Sigma_1],\dots,[\Sigma_d])$ is in $\Z$.

We consider an $SO(3)$-bundle $P$ over $X=Y \# S^2 \times S^2$ satisfying
\begin{equation*} \label{eq w_2 P}
w_2(P)=w_2(X), \quad p_1(P)=p_1(Q)-4,
\end{equation*}
so that $P$ satisfies (\ref{eq w_2}). The dimension of $M_P$ is given by $2d+5$. 

We define surfaces $\Sigma$, $\Sigma'$ embedded in $S^2 \times S^2$ by
\[
\Sigma=S^2 \times \{ pt \}, \quad \Sigma'=\{ pt \} \times S^2 \subset S^2 \times S^2.
\]
Then we have
\[
[\Sigma] \cdot [\Sigma] \equiv [\Sigma'] \cdot [\Sigma'] \equiv 0 \mod 2.
\]
Now $q_{k,w,Y \# S^2 \times S^2}^{u_1}([\Sigma_1],\dots,[\Sigma_d],[\Sigma],[\Sigma'])$ is defined for homology classes $[\Sigma_i]$ of $Y$ with self-intersection numbers even.
The following is an $SO(3)$-version of Theorem 1.1 in \cite{FS}.

\begin{theorem} \label{main thm}
In the above situation, we have
\[
q_{k,w,Y \# S^2 \times S^2}^{u_1}([\Sigma_1],\dots,[\Sigma_d],[\Sigma],[\Sigma']) \equiv q_{k-1,w,Y}([\Sigma_1],\dots,[\Sigma_d]) \mod 2.
\]
\end{theorem}

The proof is given in the following three subsections. 

\subsection{Notations and general facts} \label{facts}
For the proof of Theorem \ref{main thm}, we will investigate the intersection $M_P \cap V_{\Sigma_1} \cap \cdots \cap V_{\Sigma_d} \cap V_{\Sigma'} \cap V_{\Sigma}$ when the neck of $Y \# S^2 \times S^2$ is very long. For the preparation, we define some notations and recall some facts about instantons over the connected sum of $4$-manifolds.

Let $Y_1$ and $Y_2$ be a closed, oriented $4$-manifold.
The connected sum $X=Y_1 \# Y_2$ is constructed in the following way.
Fix Riemannian metrics $g_1$ and $g_2$ on $Y_1$ and $Y_2$ which are flat in small neighborhoods of fixed points $y_1 \in Y_1$ and $y_2 \in Y_2$. 
For $N>1$ and $\lambda>0$ with $N\lambda^{\frac{1}{2}} \ll 1$, we put
\[
\Omega_i=\Omega_{y_i}(\lambda,N)= \{ y \in Y_i | N^{-1} \lambda^{\frac{1}{2}} < d(y,y_i) < N \lambda^{\frac{1}{2}} \} \quad (i=1, 2).
\]
Let
\[
\sigma:(TY_1)_{y_1} \stackrel{\cong}{\longrightarrow} (TY_2)_{y_2}.
\]
be an orientation-reversing linear isometry. For each positive real number $\lambda>0$, we define
\[
\begin{array}{rccc}
f_{\lambda}: & (TY_1)_{y_1} \backslash \{ 0 \} & \longrightarrow & (TY_2)_{y_2} \backslash \{ 0 \} \\
& \xi & \longmapsto & \frac{\lambda}{| \xi |^2}\sigma(\xi).
\end{array}
\]
This map $f_{\lambda}$ induces a diffeomorphism between $\Omega_1$ and $\Omega_2$. The connected sum $X$ of $Y_1$ and $Y_2$ is identified with
\[
X(\lambda)=(Y_1 \backslash B_{y_1}(N^{-1} \lambda^{\frac{1}{2}}) ) \bigcup_{f_{\lambda}} (Y_2 \backslash B_{y_2}(N^{-1} \lambda^{\frac{1}{2}}))
\]
where $B_{y_i}(N^{-1} \lambda^{\frac{1}{2}})$ is the open ball centered on $y_i$ with radius $N^{-1} \lambda^{\frac{1}{2}}$. The metrics $g_1$ and $g_2$ define a conformal structure on $X$ since $g_i$ is flat in a small neighborhood of $y_i$. We fix a metric $g_{\lambda}$ on $X$ which represents the conformal structure. Moreover we assume that $g_{\lambda}$ is equal to $g_i$ on $Y_i \backslash B((N+1) \lambda^{\frac{1}{2}})$.

\begin{definition}
Fix a real number $q$ with $q>4$.
Let $[A^{(n)}] \in M_P(g_{\lambda_n})$ be instantons over $X=Y_1 \# Y_2$ for a sequence $\lambda_n \rightarrow 0$. Let $z_1,\dots, z_l$ be points in $Y_1 \backslash \{ y_1 \}$, $z_1',\dots, z_m'$ be points in $Y_2 \backslash \{ y_2 \}$ and $A_i$ be connections over $Y_i$. Then we say that $[A^{(n)}]$ is weakly convergent to $([A_1], [A_2];z_1,\dots,z_l, z_1',\dots, z_m')$ when $[A^{(n)}]$ is $L^q$-convergent to $([A_1],[A_n])$ over compact subsets in $(Y_1 \cup Y_2) \backslash \{ y_1, y_2, z_1, \dots,z_l, z_1',\dots, z_m'\}$ and $|F_{A^{(n)}}|^2$ is convergent as measure to
\[
|F_{A_1}|^2 + |F_{A_2}|^2 + 8 \pi^2 \left( \sum_{\nu=1}^l \delta_{z_{\nu}} + \sum_{\nu=1}^{m} \delta_{z_{\nu}'} \right)
\]
over compact subsets in $(Y_1 \backslash \{ y_1 \}) \cup (Y_2 \backslash \{ y_2 \})$. Here $\delta_z$ is the delta function supported on $z$.
\end{definition}

We use the following well-known theorem.

\begin{theorem} [\cite{poly, DK}] \label{thm weak conv}
Let $P$ be an $SO(3)$-bundle over $X=Y_1 \# Y_2$. Set $k=-p_1(P)/4$, $w=w_2(P)$, $w_i=w|_{Y_i}$.
Let $[A^{(n)}] \in M_{k, w, X}(\lambda_n)$ be instantons over $X$ for $\lambda_n \rightarrow 0$. Then there is a subsequence $\{ [A^{(n')}] \}_{n'}$ which is weakly convergent to $([A_1],[A_2];z_1,\dots, z_l, z_1',\dots, z_m')$ for some
\[
[A_1] \in M_{k_1,w_1,Y_1}(g_1), \ 
[A_2] \in M_{k_2,w_2,Y_2}(g_2), \
z_1,\dots, z_l \in Y_1 \backslash \{ y_1 \}, \
z_1',\dots, z_m' \in Y_2 \backslash \{ y_2 \}
\]
with
\[
k_1 \geq 0, \quad k_2 \geq 0, \quad
k_1 + k_2 + l + m \leq k.
\]
\end{theorem}

Next we review gluing of instantons. The theory of gluing of instantons is standard. To fix notations, we recall the theory briefly.

Let $A_i$ be instantons over $Y_i$. We denote the $SO(3)$-bundles carrying $A_i$ by $P_i$. We can construct instantons on $X=Y_1 \# Y_2$ close to $A_i$ on each factor. Outline of the construction is as follows. (See \cite{DK} Chapter 7 for details.)

Let $b$ be a small positive number with $b \geq 4 N \lambda^{\frac{1}{2}}$. By using suitable cut-off functions and trivializations of $P_i$ on neighborhoods of $y_i$, we obtain a connections $A_i'$ which are flat over the annuli $\Omega_i$ and equal to $A_i$ outside the balls centered at $y_i$ with radius $b$. Take an $SO(3)$-isomorphism $\rho$ between $(P_1)_{y_1}$ and $(P_2)_{y_2}$. We can spread this isomorphism by using flat structures of $A_i'$, and obtain an isomorphism $g_{\rho}$ between $P_1|_{\Omega_1}$ and $P_2|_{\Omega_2}$ covering $f_{\lambda}$. We define an $SO(3)$-bundle $P_{\rho}$ over $X$ and a connection $A'(\rho)=A_1' \#_{\rho} A_2'$ on $P_{\rho}$ by gluing $P_i$, $A_i$ through $g_{\rho}$. Then in large region outside the neck of $X$, $A'(\rho)$ satisfies the instanton equation, and $F^+_{A'(\rho)}$ is very small near the neck. To obtain a genuine instanton we have to perturb $A'(\rho)$. We consider the equation
\begin{equation} \label{eq perturb}
F_{A'(\rho) + a }^+ = 0
\end{equation}
for $a \in \Omega_X^1(\fg_{P_{\rho}})$.
To solve this equation, we take linear maps
\[
\sigma_i:H_{A_i}^2 \longrightarrow \Omega_{Y_i}^+(\fg_{P_i})
\]
such that $d_{A_i}^+ \oplus \sigma_i$ are surjective and for each $h_i \in H^2_{A_i}$ the supports of $\sigma_i(h_i)$ are in the complement of the ball centered at $y_i$ with radius $b$. Then put
\[
\sigma:=\sigma_1+\sigma_2:H^2_{A_1} \oplus H^2_{A_2} \longrightarrow \Omega_{X}^+(\fg_{P_{\rho}}).
\]
We can construct a right inverse of $d_{A'(\rho)}^+ + \sigma$ starting from right inverses of $d_{A_i}^+ + \sigma_i$ . Decompose the right inverse as $P \oplus \pi$, where
\[
P:\Omega_{X}^+(\fg_{P_{\rho}}) \longrightarrow \Omega^1(\fg_{P_{\rho}}), \quad
\pi:\Omega_X^+(\fg_{P_{\rho}}) \longrightarrow H_{A_1}^2 \oplus H_{A_2}^2.
\]
Instead of (\ref{eq perturb}), we first consider the equation 
\[
F_{A'(\rho)+a}^+ + \sigma(h)=0
\]
for $(a, h) \in \Omega_X^1(\fg_{P_{\rho}}) \times (H_{A_1}^2 \oplus H_{A_2}^2)$. We find a solution of this equation in the form $a=P \xi$, $h=-\pi \xi$. In this case, we see that the equation is equivalent to the equation
\[
\xi + (P \xi \wedge P \xi)^+ = - F_{A'(\rho)}^+
\]
by a short calculation.
Using the contraction mapping principle, we can show that there is a unique small solution $\xi_{\rho} \in \Omega^+(\fg_{P_{\rho}})$ for the equation. We get a genuine instanton if and only if $\pi \xi_{\rho}=0$. Therefore there is a map
\[
\Psi:Gl_{y_1,y_2} \longrightarrow H^2_{A_1} \times H_{A_2}^2
\]
such that the solutions of $\Psi=0$ represent instantons over $X$. Here $Gl_{y_1,y_2}$ is the space of $SO(3)$-equivariant isomorphisms between $(P_1)_{y_1}$ and $(P_2)_{y_2}$. We fix an element $\rho_0 \in Gl_{y_1, y_2}$ to identify $Gl_{y_1,y_2}$ with $SO(3)$.

We can include the deformations of $[A_i]$ to this construction. For small neighborhoods $U_{A_i}$ of $0$ in $H_{A_i}^1$, we have a map
\[
\Psi:T:=U_{A_1} \times U_{A_2} \times SO(3) \longrightarrow H_{A_1}^2 \times H_{A_2}^2
\]
such that elements of $\Psi^{-1}(0)$ correspond to instantons.

Let $\Gamma_{A_i}$ be the isotropy group of $A_i$ in the gauge group and put $\Gamma=\Gamma_{A_1} \times \Gamma_{A_2}$. We assume that $U_{A_i}$ is $\Gamma_{A_i}$-invariant. Then there are natural actions of $\Gamma$ on $T$ and on $H_{A_1}^2 \times H_{A_2}^2$. We can show that $\Psi$ is $\Gamma$-equivariant and instantons corresponding to elements of $\Psi^{-1}(0)$ are gauge equivalent to each other if and only if they are in the same $\Gamma$-orbit. Hence we can regard $\Psi^{-1}(0)/\Gamma$ as a subspace of $M_P$.

An important feature is that instantons over $X=Y_1 \# Y_2$ which is close to $A_i$ over $Y_i$ are given in the above description. More precise statement is the following:

Let $Y_i''$ be the complement of balls centered at $y_i$ with radius $\lambda^{\frac{1}{2}}/2$. Take instantons $A_i$ over $Y_i$ and a positive number $\nu>0$. Then put
\begin{equation} \label{eq U}
U_{\lambda}(\nu):= \{ \ [A] \in \cB_X^* \ | \ d_{ q}([A|_{Y_i''}], [A_i |_{Y_i''}]) < \nu, \ i=1, 2 \ \}.
\end{equation}
Here $q$ is the fixed real number with $q>4$ and $d_q$ is the distance induced by $L^q$-norm over $Y_i''$. If $\nu>0$ is small, then there is a positive number $\lambda(\nu)>0$ such that for $\lambda< \lambda(\nu)$ we can take a neighborhood $T$ of $\{ 0 \} \times \{ 0 \} \times SO(3)$ in $H_{A_1}^1 \times H_{A_2}^1 \times SO(3)$ such that $M_P(g_{\lambda}) \cap U_{\lambda}(\nu)$ is homeomorphic to $\Psi^{-1}(0)/\Gamma$.  Summing up these:

\begin{theorem} \label{thm gluing}
Let $A_1, A_2$ be instantons on $Y_1, Y_2$. Then there is a $\Gamma=\Gamma_{A_1} \times \Gamma_{A_2}$-invariant neighborhood $T$ of $SO(3) \times \{ 0 \} \times \{ 0 \}$ in $SO(3) \times H^1_{A_1} \times H_{A_2}^1$ and $\Gamma$-equivariant map
\[
\Psi:T \longrightarrow H^2_{A_1} \times H_{A_2}^2
\]
such that $\Psi^{-1}(0)/\Gamma$ is homeomorphic to an open set $N$ in $M_{P}$. 
Moreover for a small positive number $\nu>0$, there is a $\lambda(\nu)>0$ and $T$ such that if $\lambda<\lambda(\nu)$ then $N=M_P(g_{\lambda}) \cap U_{\lambda}(\nu)$.
\end{theorem}

In particular, when $Y_2$ is $S^4$ and $A_2$ is the fundamental instanton $J$ with instanton number one, we have:

\begin{corollary} \label{coro gluing}
Let $A_1$ be an instanton over $Y_1$ and $A_2$ be the fundamental instanton $J$ over $S^4$.
For a small positive number $\nu>0$,  there is a positive number $\lambda_0>0$ and a neighborhood $U_{A_1}$ of $0$ in $H^1_{A_1}$, a neighborhood $U_{0}$ of $0$ in $S^4=\R^4 \cup \{ \infty \}$ and $\Gamma=\Gamma_{A_1}$-equivariant map
\[
\Psi:U_{A_1} \times U_{0} \times (0, \lambda_0) \times SO(3) \longrightarrow H_{A_1}^2
\]
such that $\Psi^{-1}(0)/\Gamma$ is naturally homeomorphic to $M_P \cap U_{\lambda_0}(\nu)$.
\end{corollary}

\begin{remark} \label{rem gluing}
We can generalize the statements of Theorem \ref{thm gluing} and Corollary \ref{coro gluing} to $3$ or more instantons.
\end{remark}

\subsection{Shrinking the neck} \label{neck}
In the situation of Theorem \ref{main thm}, we investigate 
\[
M_P(g_{\lambda}) \cap V_{\Sigma_1} \cap \cdots \cap V_{\Sigma_d} \cap V_{\Sigma} \cap V_{\Sigma'}
\]
as $\lambda$ tends to $0$. We use the notations in \S \ref{facts}.

Let $Y_1$ be a closed, oriented, simply connected, non-spin $4$-manifold with $b^+(Y_1)=2a-1$ with $a>1$ and we write $Y_2$ for $S^2 \times S^2$. Let $P$ be an $SO(3)$-bundle over $X=Y_1 \# Y_2$ satisfying (\ref{eq w_2}). Assume that the virtual dimension of $M_P$ is $2d+5$ for a non-negative integer $d$. Take homology classes $[\Sigma_1], \dots, [\Sigma_d] \in H_2(Y_1; \Z)$ with $[\Sigma_i] \cdot [\Sigma_i] \equiv 0 \mod 2$. Set $\Sigma=S^2 \times \{ pt \}, \Sigma'=\{ pt \} \times S^2 \subset Y_2$.
Take instantons
\[
[A^{(n)}] \in M_P(g_{\lambda_n}) \cap V_{\Sigma_1} \cap \cdots \cap V_{\Sigma_d} \cap V_{\Sigma} \cap V_{\Sigma'}
\]
for a sequence $\lambda_n \rightarrow 0$. By Theorem \ref{thm weak conv}, a subsequence of $\{ [A^{(n)}] \}_{n}$ is weakly convergent to some
\[
([A_1], [A_2]; z_1,\dots, z_{l}, z_1',\dots, z_{m}'),
\]
where 
\[
[A_1] \in M_{k_1,w,Y_1}(g_1), \ [A_2] \in M_{k_2,Y_2}(g_2), \
z_1,\dots, z_l \in Y_1 \backslash \{ y_1 \}, \
z_1',\dots , z_m' \in Y_2 \backslash \{ y_2 \}.
\]

\begin{lemma} \label{lem bubble}
In the above situation, we have
\begin{gather*}
k_1=k-1,\ l=0, \ [A_1] \in M_{k-1,w,Y_1}(g_1) \cap V_{\Sigma_1} \cap \cdots \cap V_{\Sigma_d}, \\
m=1, \quad  z_1' \in \nu(\Sigma) \cap \nu(\Sigma'), \quad [A_2]=[\Theta_{Y_2}].
\end{gather*}
Here $\Theta_{Y_2}$ is the trivial connection on $Y_2$.
\end{lemma}

\begin{proof}
From Theorem \ref{thm weak conv}, we have
\begin{equation} \label{eq p1}
k_1+k_2+l+m \leq k.
\end{equation}
Let $p$ be the number of $\nu(\Sigma_i)$ which contain some point $z_{\alpha}$ and $q$ be the number of $\nu(\Sigma)$, $\nu(\Sigma')$ which contain some point $z_{\alpha}'$. Then by the transversality condition (\ref{eq transverse}), we have
\begin{equation} \label{eq p q}
0 \leq p \leq 2l,  \quad 0 \leq q \leq 2m.
\end{equation}
Without loss of generality, we may assume
\[
[A_1] \in M_{k_1,w,Y_1} \cap V_{\Sigma_1} \cap \cdots \cap V_{\Sigma_{d-p}}
\]
if we change the order of surfaces.
Since $w_2(P)|_{Y_1}$ is non-trivial, we can show $k_1>0$ in the same way as the proof of Lemma \ref{lem compact}. For generic sections, the intersection $M_{k_1,w,Y_1} \cap V_{\Sigma_1} \cap \cdots \cap V_{\Sigma_{d-p}}$ is transverse by Lemma \ref{lem transverse}. Hence we have
\begin{equation} \label{eq M 1}
2(d-p) \leq \dim M_{k_1,w,Y_1}.
\end{equation}
We would like to show $k_2=0$. Suppose that $k_2$ is positive. Then we also obtain
\begin{equation} \label{eq M 2}
2(2-q) \leq \dim M_{k_2,Y_2}.
\end{equation}
By index theorem, there is the formula
\begin{equation} \label{eq sum}
\dim M_{k_1,w,Y_1} + \dim M_{k_2,Y_2}+3=\dim M_{k_1+k_2,w,X}.
\end{equation}
From (\ref{eq p1}), (\ref{eq M 1}), (\ref{eq M 2}) and (\ref{eq sum}), we have
\[
2(d-p)+2(2-q)+3 \leq \dim M_{k_1+k_2,w,X} \leq \dim M_{k,w,X}-8(l+m)=2d+5-8(l+m).
\]
This inequality and (\ref{eq p q}) imply
\[
8(l+m) +2 \leq 2p+2q \leq 4(l+m).
\]
We have a contradiction. Hence $k_2$ is $0$ which implies that $[A_2]$ is the class of trivial flat connection $[\Theta_{Y_2}]$.

Since $k_2$ is $0$, the virtual dimension of $M_{0,Y_2}$ is $-6$. From (\ref{eq sum}), we have
\begin{equation} \label{eq sum 2}
\dim M_{k_1,w,Y_1}-3 = \dim M_{k_1,w,X}.
\end{equation}
By (\ref{eq p1}), (\ref{eq p q}),(\ref{eq M 1}) and (\ref{eq sum 2}), we have
\[
2(d-2l) - 3 \leq 2(d-p) - 3 \leq \dim M_{k_1,w,Y_1}-3= \dim M_{k_1,w,X} \leq \dim M_{k,w,X} - 8(l+m).
\]
Therefore we obtain
\[
4l+8m \leq 8.
\]
In particular, we have $m \leq 1$. We show $m=1$. Suppose $m=0$, then we have $[\Theta_{Y_2}] \in V_{\Sigma}, [\Theta_{Y_2}] \in  V_{\Sigma'}$. To obtain a contradiction, we need to choose $V_{\Sigma}$ and $V_{\Sigma'}$ in a specific way. As mentioned in Remark \ref{rem line bundle}, we can choose $V_{\Sigma}$ and $V_{\Sigma'}$ do not include $[\Theta_{Y_2}]$. If we choose such $V_{\Sigma}$ and $V_{\Sigma'}$, we have a contradiction. We obtain $l=0$, $m=1$ and $z_1' \in \nu(\Sigma) \cap \nu(\Sigma')$. Hence 
\[
[A_1] \in M_{k_1,w,Y_1} \cap V_{\Sigma_1} \cap \cdots \cap V_{\Sigma_d}.
\]
Lastly we show $k_1=k-1$. From (\ref{eq p1}), we have $k_1 \leq k-1$. On the other hand, from (\ref{eq M 1}) we have
\[
2d \leq \dim M_{k_1,w,Y_1}=\dim M_{k-1,w,Y_1}-8(k-1-k_1)=2d-8(k-1-k_1).
\]
This implies $k_1 \geq k-1$.  Therefore $k_1$ is equal to $k-1$. We complete the proof.
\end{proof}

Let $w_0'$ be the unique intersection point of $\Sigma$ and $\Sigma'$. Fix a small neighborhood $U_{w_0'}$ of $w_0'$ with $\nu(\Sigma) \cap \nu(\Sigma') \subset U_{w_0'}$.
We suppose that the metric $g_2$ on $Y_2$ is flat on $U_{w_0'}$ for simplicity.

Take 
\[
[A^{(n)}] \in M_P(g_{\lambda_n}) \cap V_{\Sigma_1} \cap \cdots \cap V_{\Sigma_d} \cap V_{\Sigma} \cap V_{\Sigma'}
\]
for $\lambda_n \rightarrow 0$ and assume that $\{ [A^{(n)}] \}_{n \in \N }$ weakly converges to $([A_1], [\Theta_{Y_2}]; z_1')$ for some $[A_1] \in M_{k-1,w,Y_1} \cap V_{\Sigma_1} \cap \cdots \cap V_{\Sigma_d}$,  $z_1' \in \nu(\Sigma) \cap \nu(\Sigma')$. We can define the local center of mass $c_n \in U_{w_0'}$ and scale $\lambda_n'>0$ of $[A^{(n)}]$ around $z_1'$ when $n$ is sufficiently large.  If $n$ is large enough, then we obtain
\[
\int_{U_{w_0'}} | F_{A^{(n)}}|^2 d\mu_{g_2} > 4\pi^2
\]
since $|F_{A^{(n)}}|^2$ converges to $8\pi^2 \delta_{z_1'}$ on $U_{w_0'}$.
We define the center of mass $c_n$ to be the center of the smallest ball in $U_{w_0'}$ where the integral of $|F_{A^{(n)}}|^2$ is equal to $4 \pi^2$ and the scale $\lambda_n'$ to be the radius of the ball. The center of mass and scale is determined uniquely (\cite{appli}). The center $c_n$ converges to $z_1'$ and the scale $\lambda_n'$ converges to $0$. 

Let $m:\R^4 \rightarrow S^4=\R^4 \cup \{ \infty \}$ be the stereographic map and $d_{\lambda}:\R^4 \rightarrow \R^4$ be the map $d_{\lambda}(y)=\lambda^{-1}y$. Put $\chi_n:=m \circ d_{\lambda_n'}$. 
Then $\chi_n$ induces a conformal isomorphism between $X$ and the connected sum \[
X \# S^4=(X \backslash B_{c_n}(N^{-1}\lambda_n')) \cup_{f_{\lambda_n'}} (S^4 \backslash B_{\infty}(N^{-1} \lambda_n') )
\]
since the metric $g_2$ is flat on $U_{w_0'}$. Here $f_{\lambda_n'}$ is defined in the following way: Using the geodesic coordinate near $c_n$ and the stereographic map, we identify $(TX)_{c_n}$ with $(TS^4)_{0}$. Let $\sigma'$ be the natural, orientation reversing isometry between $(TS^4)_{0}$ and $(TS^4)_{\infty}$, then $f_{\lambda_n'}$ is given by
\[
\begin{array}{rccc}
f_{\lambda_n'}: & (TX)_{c_n} \backslash \{ 0 \} & \longrightarrow & (TS^4)_{\infty} \backslash \{ 0 \} \\
& \xi & \longmapsto & \frac{\lambda_n'}{|\xi|^2} \sigma'(\xi).
\end{array}
\]
We can regard $A^{(n)}$ as an instanton on $X \# S^4$ such that $A^{(n)}$ is close to $A_1$, $\Theta_{Y_2}$ on $Y_1$, $Y_2$ and close to the standard instanton  $J$ on $S^4$. 

Fix a small positive number $\lambda_0$ and a small neighborhood $U_{[A_1]}'$ of $[A_1]$ in $M_Q$.
Let $O_{[A_1]} \subset \cB_{P}^*$ be a small open neighborhood of
\[
\{ \ [B' \ \#_{y_1,\lambda, \rho} \ \Theta_{Y_2} \ \#_{z_1', \lambda', \rho'} \ J' ] \ | \
B \in U_{[A_1]}', \ \lambda, \lambda' \in (0, \lambda_0), \ \rho, \rho' \in SO(3), \ z_1' \in \nu(\Sigma) \cap \nu(\Sigma') \ \}.
\]
Here $B', J'$ are connections which are flat near $y_1, \infty$ and equal to $B, J$ outside $b$-balls. (The real number $b$ is a small positive number fixed in \S \ref{facts}). The instanton $[A^{(n)}]$ is in $O_{[A_1]}$ when $n$ is large.
We can define the local centers for elements of $O_{[A_1]}$ and we have a map $p:O_{[A_1]} \rightarrow U_{w_0'}$ which maps connections to their centers.
By Donaldson \cite{conn}  Proposition (3.18), we can take sections $s_{\Sigma}$, $s_{\Sigma'}$ such that $O_{[A_1]} \cap V_{\Sigma}$, $O_{[A_1]} \cap V_{\Sigma'}$ are equal to $p^{-1}(U_{z_1'}' \cap \Sigma)$, $p^{-1}( U_{z_1'}' \cap \Sigma')$. Hence we may suppose that the center $c_n$ of $[A^{(n)}]$ is $w_0'$ for large $n$.

We denote $S^4$ by $Y_3$ and denote $\Theta_{Y_2}$, $J$ by $A_2$, $A_3$ and put
\[
Y_{1,n}''=Y_1 \backslash B_{y_1}(\lambda_n / 2), \quad
Y_{2,n}''=Y_2 \backslash ( B_{y_2}(\lambda_n / 2) \cup B_{w_0'}(\lambda_n' / 2)), \quad
Y_{3,n}''=Y_3 \backslash B_{\infty}(\lambda_n' / 2).
\]
For $\nu>0$, put
\[
U_{[A_1],\lambda_n}(\nu)=
\{ \ [A] \in \cB_{X \# S^4}^* \ | \ d_{q}([A|_{Y_{i,n}''}], [A_i|_{Y_{i,n}''}])< \nu, \ i=1, 2, 3 \ \}.
\]
We have proved the following:

\begin{lemma}
Fix a positive number $\nu>0$.
Take instantons $[A^{(n)}] \in M_P(g_{\lambda_n}) \cap V_{\Sigma_1} \cap \cdots \cap V_{\Sigma_d} \cap V_{\Sigma} \cap V_{\Sigma'}$ for a sequence $\lambda_n \rightarrow 0$. Then $[A^{(n)}]$ is in $U_{[A_1], \lambda_n}(\nu)$ for some $[A_1] \in M_Q \cap V_{\Sigma_1} \cap \cdots \cap V_{\Sigma_d}$ when $n$ is sufficiently large.
\end{lemma}

Fix $[A_1] \in M_{Q} \cap V_{\Sigma_1} \cap \cdots \cap V_{\Sigma_d}$ and a small positive number $\nu$. By Theorem \ref{thm gluing}, Corollary \ref{coro gluing} and Remark \ref{rem gluing},  there is a small neighborhood $U_{A_1}$ of $0$ in $H_{A_1}^1$, a positive real number $\lambda_0$ and a $\Gamma_{\Theta_{Y_2}}$-equivariant map
\[
\Psi: T=U_{A_1} \times SO(3) \times U_{w_0'} \times (0, \lambda_0) \times SO(3) \longrightarrow H^2_{\Theta_{Y_2}}
\]
such that $\Psi^{-1}(0)/\Gamma_{\Theta_{Y_2}}$ is homeomorphic to $M_P(g_{\lambda_n}) \cap U_{[A_1], \lambda_n}(\nu)$. Since the action of $\Gamma_{\Theta_{Y_2}}=SO(3)$ on $SO(3) \times SO(3)$ is the diagonal action, $\Psi^{-1}(0)/SO(3)$ is naturally identified with
\[
\Psi^{-1}(0) \cap \big( U_{A_1} \times \{ 1 \} \times U_{w_0'} \times (0, \lambda_0) \times SO(3) \big).
\]
We write $T'$ for $U_{A_1} \times \{ 1 \} \times U_{w_0'} \times (0, \lambda_0) \times SO(3)$.
Since $T'$ parametrizes connections on $X$, it makes sense to take the intersection $T' \cap V_{\Sigma_1} \cap \cdots \cap V_{\Sigma_d} \cap V_{\Sigma} \cap V_{\Sigma'}$. We can suppose 
\[
T' \cap V_{\Sigma_1} \cap \cdots \cap V_{\Sigma_d} \cap V_{\Sigma} \cap V_{\Sigma'} =
\{ 0 \} \times \{ 1 \} \times \{ w_0' \} \times (0,\lambda_0) \times SO(3).
\]
Hence $M_P(g_{\lambda_n}) \cap V_{\Sigma_1} \cap \cdots \cap V_{\Sigma_d} \cap V_{\Sigma} \cap V_{\Sigma'} \cap U_{[A_1],\lambda_n}(\nu)$ is homeomorphic to
\[
\Psi^{-1}(0) \cap \big( \{ 0 \}  \times \{ 1 \} \times \{ w_0' \} \times (0,\lambda_0) \times SO(3) \big)
\subset 
H^1_{A_1} \times SO(3) \times U_{w_0'} \times (0,\lambda_0) \times SO(3).
\]
Donaldson calculated the leading term of $\Psi$ in \cite{conn} explicitly. By the explicit expression of the leading term of $\Psi$ and calculations similar to those in \cite{conn} V, we can show the following:

\begin{lemma}
For generic metrics $g_1$ and $g_2$, and points $y_1$, $y_2$ and $w_0'$, the intersection 
\[
\Psi^{-1}(0) \cap \big( \{ 0 \}  \times \{ 1 \} \times \{ w_0' \} \times (0,\lambda_0) \times SO(3) \big)
\]
is homeomorphic to
\[
\{ c \lambda_{n} \} \times \gamma \subset (0, \lambda_0) \times SO(3)
\]
where $\gamma$ is a loop in $SO(3)$ which represent the generator of $\pi_1(SO(3)) \cong \Z_2$ and $c>0$ is a constant number independent of $n$.
\end{lemma}

Define $N_{[A_1]}$ by 
\begin{equation} \label{eq NA}
N_{[A_1]}=\{ \ [A_1' \#_{\lambda_n} \Theta_{Y_2} \#_{w_0', c \lambda_n, \rho} J'] \ | \ \rho \in \gamma \ \}.
\end{equation}
We have obtained the following:

\begin{corollary} \label{coro intersection}
Let $Y$ be a closed, oriented, simply connected, non-spin $4$-manifold with $b^+(Y)=2a-1$ for $a>1$ and $P$ be an $SO(3)$-bundle over $X=Y \# S^2 \times S^2$ which satisfies the condition (\ref{eq w_2}). Suppose that the virtual dimension of $M_P$ is $2d+5$ for a non-negative integer $d$. Take $d$ homology classes $[\Sigma_i]$ in $H_2(Y; \Z)$ with self-intersection numbers even. Then for a small $\lambda>0$, generic metrics $g_1$ and $g_2$, and generic points $y_1,y_2$ and $w_0'$, the intersection
\[
M_P(g_{\lambda}) \cap V_{\Sigma_1} \cap \cdots \cap V_{\Sigma_d} \cap V_{\Sigma} \cap V_{\Sigma'}
\]
is homeomorphic to
\[
\coprod_{[A_1] \in M_Q \cap V_{\Sigma_1} \cap \cdots \cap V_{\Sigma_d}} N_{[A_1]}.
\]
\end{corollary}

\subsection{End of the proof}
From Corollary \ref{coro intersection}, we have
\[
q_{k,w,Y \# S^2 \times S^2}^{u_1}([\Sigma_1],\dots,[\Sigma_d], [\Sigma], [\Sigma'])= \sum_{[A_1]} \< u_1, N_{[A_1]} \> \in \Z_2,
\]
where $[A_1]$ runs in $M_Q \cap V_{\Sigma_1} \cap \cdots \cap V_{\Sigma_d}$. Therefore it is sufficient to show that the pairing $\< u_1, N_{[A_1]} \>$ is non-trivial for the proof of Theorem \ref{main thm}. The last step is carried out by making use of the following Proposition due to Akbulut, Mrowka and Ruan.

\begin{proposition} [\cite{AMR}] \label{prop u1}
Let $X_i$ be closed, oriented, simply connected $4$-manifolds for $i=1,2$ and $x_i$ be points of $X_i$. Take $SO(3)$-bundles $P_i$ over $X_i$ with $w_2(P_i)$ equal to $w_2(X_i)$. Choose $U(2)$-lifts $\bar{P}_i$ of $P_i$ and  assume that the second Chern numbers of $\bar{P}_i$ are odd. (In this case, $P_1 \# P_2$ satisfies the condition (\ref{eq w_2}). See  Remark \ref{rem pi1}.) We fix trivializations of $P_i$ on small neighborhoods $U_{x_i}$ of $x_i$. For irreducible connections $B_i$ on $P_i$ with trivial on $U_{x_i}$ with respect to fixed trivializations, we have a family of connections
\[
G:=\{ \ [B_1 \#_{\rho} B_2] \ | \ \rho \in SO(3) \ \} \ (\cong SO(3)) 
\subset \cB_{P_1 \# P_2}^*.
\] 
 Then the restriction $u_1|_G$ is non-trivial in $H^1(G; \Z_2) \cong \Z_2$.
\end{proposition}

In our case, 
\[
X_1=Y \# S^2 \times S^2,\ P_1=Q \# P_{S^2 \times S^2}, \ B_1=A_1' \# \Theta_{S^2 \times S^2},  \
X_2=S^4, \ P_2=P_{S^4}/\{ {\pm 1} \}, \ B_2=J'.
\]
Here $Q$ is an $SO(3)$-bundle over $Y$ with
\begin{equation} \label{eq Q w2}
w_2(Q)=w_2(Y), \quad p_1(Q) \equiv \sigma (Y)+4 \mod 8,
\end{equation}
$P_{S^2 \times S^2}$ is the trivial $SO(3)$-bundle over $S^2 \times S^2$ and $P_{S^4}$ is an $SU(2)$-bundle with second Chern number equal to $1$. By the formulas
\[
p_1(Q)=-4c_2(\bar{Q}) + c_1(\bar{Q})^2, \ w_2(Y)^2 \equiv \sigma(Y) \mod 8
\]
and (\ref{eq Q w2}), we have
\[
c_2(\bar{Q}) \equiv 1 \mod 2.
\]
Hence the assumptions of Proposition \ref{prop u1} is satisfied. Since $N_{[A_1]}$ is a loop in $G$ which represent the generator of $\pi_1(G) \cong \Z_2$, we obtain:

\begin{corollary}
For each $[A_1] \in M_Q \cap V_{\Sigma_1} \cap \cdots \cap V_{\Sigma_d}$, the pairing $\< u_1, N_{[A_1]} \>$ is non-trivial in $\Z_2$. 
\end{corollary}

This completes the proof of Theorem \ref{main thm}.

\section{Example}

\subsection{Non-triviality of $q_{2\CP^2 \# \barCP2}^{u_1}$} \label{example}
We see that the $SO(3)$-torsion invariant for $X=2\CP^2 \# \overline{\CP}^2$ is non-trivial. 

To distinguish two $\CP^2$'s, we write $X=\CP^2_1 \# \CP^2_2 \# \barCP2$.

\begin{theorem} \label{thm non-trivial}
Let $H_i$ be the canonical generator of $H_2(\CP^2_i;\Z)$ for $i=1, 2$ and $E$ be the canonical generator of $H_2(\barCP2; \Z)$. Then we have
\[
q_{\CP^2_1 \# \CP_2^2 \# \barCP2}^{u_1}(-H_1 + E, H_2 - E) \equiv 1 \mod 2.
\]
\end{theorem}

\begin{proof}
Let $Q$ be an $SO(3)$-bundle on $\CP^2$ with
\[
w_2(Q)=w_2(\CP^2), \quad p_1(Q)=-3.
\]
Then the dimension of $M_Q$ is $0$.  Kotschick showed that the Donaldson invariant associated with $Q$ is 
\[
q_{\frac{3}{4},w,\CP^2}=-1
\]
if we choose a suitable orientation on $M_Q$ (\cite{K, K2}).  Note that there is no wall since $b^-(\CP^2)$ is $0$. 
The signature of $\CP^2$ is $1$, hence we have
\[
p_1(Q) \equiv \sigma(\CP^2) + 4 \mod 8
\]
and $q^{u_1}_{\frac{7}{4},w,\CP^2 \# S^2 \times S^2}([\Sigma],[\Sigma'])$ is defined. From Theorem \ref{main thm}, we have
\[
q^{u_1}_{\frac{7}{4},w,\CP^2 \# S^2 \times S^2}([\Sigma], [\Sigma']) \equiv 1 \mod 2.
\]
On the other hand, $\CP^2 \# S^2 \times S^2$ is diffeomorphic to $\CP^2_1 \# \CP^2_2 \# \barCP2$ (\cite{W}). The induced isomorphism between the $2$-dimensional homology groups is given by
\[
\begin{array}{ccl}
H_2(\CP^2 \# S^2 \times S^2; \Z) &\stackrel{\cong}{\longrightarrow} & H_2(\CP^2_1 \# \CP_2^2 \# \barCP2; \Z) \\
H & \longmapsto & H_1 + H_2 - E \\
\left[ \Sigma \right]  & \longmapsto & -H_1 + E \\
\left[ \Sigma' \right] & \longmapsto & H_2 - E.
\end{array}
\]
The torsion cohomology class $w$ is $w_2(\CP^2 \# S^2 \times S^2)$, and the image of $w$ under the isomorphism is $w_2(2\CP^2 \# \barCP2)$. We also denote this class by $w$. The images of $[\Sigma]$ and $[\Sigma']$ under the isomorphism are $-H_1 + E$ and $H_2 - E$ respectively. Hence we obtain
\[
q^{u_1}_{\frac{7}{4}, w, \CP^2_1 \# \CP^2_2 \# \barCP2}(-H_1 + E, H_2-E) \equiv 1 \mod 2.
\]
\end{proof}

\subsection{A vanishing theorem}
Let $X$ be a closed, oriented, simply connected, non-spin $4$-manifold with $b^+(X)=2a$ for some $a>0$. Moreover assume that $X$ can be written as the connected sum $Y_1 \# Y_2$ of non-spin $4$-manifolds $Y_i$ with $b^+(Y_i) \geq 1$. In this situation, we can show a vanishing theorem similar to the usual Donaldson invariant. However we must require a certain condition for homology classes in $X$. The condition is that each homology class lies in $H_2(Y_1; \Z)$ or $H_2(Y_2; \Z)$.

Suppose that $P$ is an $SO(3)$-bundle over $X$ satisfying (\ref{eq w_2}) and that $\dim M_P$ is $2d+1$ for some non-negative integer $d$. Moreover suppose that $d=d_1+d_2$ for some $d_1 \geq 0$, $d_2 \geq 0$. Take homology classes $[\Sigma_1], \dots, [\Sigma_{d_1}] \in H_2(Y_1; \Z)$, $[\Sigma_1'], \dots, [\Sigma_{d_2}'] \in H_2(Y_2; \Z)$ with self-intersection numbers even.  Then by the standard dimension-count argument \cite{MM}, we can show 
\[
M_P \cap V_{\Sigma_1} \cap \cdots \cap V_{\Sigma_{d_1}} \cap V_{\Sigma_1'} \cap \cdots \cap V_{\Sigma_{d_2}'} = \emptyset
\]
when the neck is sufficiently long. Hence we have:

\begin{theorem} \label{thm vanishing}
Let $Y_1, Y_2$ be closed, oriented, simply connected, non-spin $4$-manifolds with $b^+(Y_i)>0$ and $b^+(Y_1) \equiv b^+(Y_2) \mod 2$.  Then for homology classes $[\Sigma_1], \dots, [\Sigma_{d_1}] \in H_2(Y_1; \Z)$, $[\Sigma_1'], \dots, [\Sigma_{d_2}'] \in H_2(Y_2; \Z)$ with self-intersection numbers even, we have
\[
q_{Y_1 \# Y_2}^{u_1}([\Sigma_1], \dots, [\Sigma_{d_1}], [\Sigma_1'], \dots, [\Sigma_{d_2}']) \equiv 0 \mod 2.
\]
\end{theorem}

\begin{remark}
We regard $X=2\CP^2 \# \overline{\CP}^2$ as the connected sum of $Y_1=\CP^2$ and $Y_2=\CP^2 \# \overline{\CP}^2$. Then $w$ is non-trivial on $Y_i$ for $i=1, 2$. By Theorem \ref{main thm}, $q^{u_1}_{Y_1 \# Y_2}(-H_1+E, H_2-E)$ is non-trivial in contrast to Theorem \ref{thm vanishing}.
If there were a formula like
\begin{gather*}
q_{\frac{7}{4},w,Y_1 \# Y_2}^{u_1}(-H_1+E,H_2-E) \equiv  \\
`` q_{\frac{7}{4},w,Y_1 \# Y_2}^{u_1}(-H_1,H_2-E)" + ``q_{\frac{7}{4},w,Y_1 \# Y_2}^{u_1}(E,H_2-E)"  \mod 2,
\end{gather*}
then we would be able to apply Theorem \ref{thm vanishing} to showing the vanishing of $q_{\frac{7}{4},w,Y_1 \# Y_2}^{u_1}(-H_1+E,H_2-E)$. However $`` q_{\frac{7}{4},w,Y_1 \# Y_2}^{u_1}(-H_1,H_2-E)"$ nor $``q_{\frac{7}{4},w,Y_1 \# Y_2}^{u_1}(E,H_2-E)"$ are not defined because
\[
(-H_1) \cdot (-H_1) \equiv E \cdot E \equiv 1 \mod 2.
\]
\end{remark}

\end{document}